\documentclass[twocolumn]{autart}    

\usepackage{amsmath,amssymb,amsfonts,bm}
\usepackage{graphicx}
\usepackage{algorithm}
\usepackage{algorithmicx,algpseudocode}
\usepackage{cite}
\usepackage{color}
\let\theoremstyle\relax
\usepackage{amsthm}
\usepackage{cuted}
\usepackage[mathscr]{euscript}
\theoremstyle{plain}
\newtheorem{theorem}{Theorem}[section]
\newtheorem{lemma}[theorem]{Lemma}

\newtheorem{proposition}[theorem]{Proposition}
\theoremstyle{definition}

\newtheorem{problem}[theorem]{Problem}
\theoremstyle{remark}
\newtheorem*{remark}{Remark}
\usepackage{wrapfig}

\begin{document}

\begin{frontmatter}

\title{Group Steering: Approaches Based on Power Moments} 
\thanks[footnoteinfo]{This work was initiated and the first version was completed while the first author was with Shanghai Jiao Tong University. The first author acknowledges support from Nanyang Technological University during the revision stage. A preliminary version of this work was presented at the IFAC World Congress 2023, Yokohama, Japan~\cite{wu2023density}.}

\author[Guangyu]{Guangyu Wu}\ead{gwu1@alumni.nd.edu},    
\author[{Anders1}]{Anders Lindquist}\ead{alq@kth.se}              

\address[Guangyu]{Department of Automation, Shanghai Jiao Tong University, China}  
\address[Anders1]{School of Artificial Intelligence, Anhui University, Hefei, China}
          
\begin{keyword}                           
Distribution steering; power moments; multiple agents.               
\end{keyword}                             

\begin{abstract}                          
This paper considers the problem of steering a vast group of agents of which the dynamics are governed by a discrete-time asymptotically stable first-order linear system. The group of agents are characterized as a probability density function and an occupation measure respectively in the paper and two corresponding treatments are given. We propose to use the power moments to characterize the density function/occupation measure of the agents. A moment system representation of the original system is put forward for control and an empirical control scheme corresponding to it is proposed. By the designed control law, the moment sequence of the control at each time step is positive, which ensures the existence of the control for the moment system. We then realize the control as an analytic form of function by a convex optimization scheme of which the existence and uniqueness of the solution have been proved in our previous paper. The terminal density is proved to converge to the desired terminal one, which distinguishes the proposed distribution steering scheme from other existing ones. An error analysis of the terminal density from the specified one is also provided. For the problem where the group of agents is characterized as an occupation measure, the control for each agent is determined by drawing independent and identically-distributed(i.i.d) samples from the realized analytic function. Finally, simulation results validate our proposed algorithms.
\end{abstract}

\end{frontmatter}

\section{Introduction}
\label{sec:introduction}

We are interested in the problem of steering the states of a vast group of agents, of which the dynamics are governed by a discrete-time asymptotically stable first-order linear system, between an initial and a final probability densities or occupation measures without stochastic disturbance. We consider the linear dynamics of the $i_{\text{th}}$ agent
\begin{equation}
x_{i}(k+1) = a(k)x_{i}(k) + u_{i}(k), \quad i = 1, \cdots, N,
\label{uniequation}
\end{equation}
where $x_{i}(k), u_{i}(k) \in \mathbb{R}$. Since the system is asymptotically stable, we have
\(\lvert a(k)\rvert<1\) (see, e.g., Section~3.5.5 of \cite{antsaklis2007linear}).  Together with the assumption \(a(k)>0\),
this implies \(0<a(k)<1\) for all \(k\) .
The control input on the $i_{\text{th}}$ agent is denoted as $u_{i}$, and $x_{i}$ is its state. We assume that the agents are non-interactive and the volume of the agents is ignored. It means that the agents are allowed to occupy the same state and the collisions are ignored.

This poses a significant challenge both in theory and practice. From a theoretical standpoint, the challenge involves controlling numerous agents, distinguishing it from conventional control problems where only one object is typically controlled, and feedback control laws are developed based on the state of that single object. However, applying such principles to a vast group of agents becomes highly complex. As the number of agents approaches infinity, the problem becomes infinite-dimensional, making it intractable without dimension reduction or approximation. This intrinsic infinite-dimensionality adds depth and interest to the problem, presenting it as an open challenge. In practical terms, there are numerous scenarios demanding the control of a group of agents. These applications extend to steering swarms (such as UAVs or large collections of microsatellites), modeling the flow and collective motion of agents, and other ensemble control situations.

Numerous studies have addressed the challenge of steering groups of agents, leading to diverse findings. One line of research formulates the issue as stabilizing a discrete-time Markov process evolving within a compact subset of $\mathbb{R}^{d}$ towards a target density \cite{biswal2020stabilization, biswal2021decentralized, elamvazhuthi2023density}. These contributions are profound and thought-provoking. However, the results primarily focus on stabilizing the system and lack the capability to steer the initial density to an arbitrary distribution in a specified limited number of steps. In prior research, another notable exploration revolves around the assumption that the distribution follows a Gaussian pattern \cite{okamoto2019optimal, balci2020covariance, pilipovsky2021covariance, liu2023optimal, bakolas2020greedy, saravanos2022distributed}. Specifically, the goal here is to manipulate the first two moments of Gaussian distributions to converge towards predetermined targets. However, these Gaussian methodologies primarily originate from the context of quantifying control uncertainty, where the distributions serve to characterize the uncertainty inherent in the system state. Nevertheless, the theory developed for the task of group steering in the paper should include Gaussian distribution as a special case, but not limited to it. Expanding the probability distribution to a broader class of functions significantly enhances the effectiveness of the proposed group steering algorithm in real-world applications. Considerable advancements in steering continuous-time systems have been made, notably by Chen, Georgiou, and Pavon, who introduced fundamental results using the Schr{\"o}dinger Bridge strategy for Gaussian distributions \cite{chen2015optimal, chen2015optimal2} and extended their findings to other types of distributions \cite{chen2018optimal}. Additionally, a robust optimal density control strategy for robotic swarms is proposed in \cite{sinigaglia2022robust}. These contributions collectively enrich our understanding of steering problems across different system dynamics and constraints.

When the density of the agents is assumed to be Gaussian, we use the first and second order moments to characterize the density, which turns the problem into a finite-dimensional one. However, by generalizing the mean and covariance to all the power moments, we will have a more conceptual view of this problem. Controlling the system state as a distribution function, if only assumed to be Lebesgue integrable, is an uncountably infinite-dimensional problem. By probability theory, we note that a distribution function can be uniquely determined by its full power moment sequence \cite{Shiryaev2016, chung2001course}. By controlling the full power moment sequence instead of the distribution of system state, the problem is reduced to a countably infinite-dimensional one. By properly truncating the first several terms of the power moment sequence for characterizing the density of the system state, the problem is now steering a truncated power moment sequence to another, which is finite-dimensional and tractable. Decomposing the distribution by the power moments provides a completely new perspective to the challenge proposed by Brockett in \cite{brockett2007optimal}.

The paper is organized as follows. In Section 2, we prove a moment system representation as a counterpart of the discrete-time linear system. Subsequently, we introduce a problem formulation for group steering based on the moment system. Unlike conventional control problems, the positivity of Hankel matrices for both control inputs and system states' moments is required, posing challenges for standard approaches such as optimal control. In the context of general distribution steering, as considered in \cite{sivaramakrishnan2022distribution}, the system trajectory is predetermined. We propose an empirical scheme for determining the system trajectory. Additionally, we employ a density parametrization algorithm from our prior work \cite{wu2023non} to realize control inputs as analytic functions using power moments. With the number of moment terms used approaching infinity, we propose that the terminal distribution obtained by the proposed steering scheme converges almost everywhere to the desired one. The density steering problem is intrinsically infinite-dimensional, where densities are not assumed to adhere to specific functions. Since the number of power moments we utilize is finite, the presence of errors in the terminal density is inevitable. Consequently, we propose a tight upper bound for the error of the terminal density relative to the specified one. Building upon the density-steering algorithm outlined in Section 2, we introduce an algorithm for steering a finite number of agents characterized as an occupation measure in Section 3. The control input for each agent is generated by sampling independent and identically distributed (i.i.d.) values from the realized control functions, facilitated by an acceptance-rejection sampling technique. It ensures the high scalability of the algorithm, even when handling a significant number of agents. In Section 4, we provide four examples to validate the effectiveness of the two proposed steering algorithms. The examples encompass complicated mixtures of distributions.

\section{Steer the group of agents as a probability density function}

Let $\mathcal{D}$ denote the space of probability density functions $p$ on the real line with support there. Let $\mathcal{D}_{2n}$ be the subset of all $p \in \mathcal{D}$ that have at least $2n$ finite moments (in addition to the zeroth power moment, which is 1 by definition). With a vast group of agents, i.e., $N$ is large, a conventional approach is to approximate $x(k)$ and $u(k)$ as random variables for which the density functions are denoted as $q_{k}$ and $p_{k}$, both in $\mathcal{D}_{2n}$. It is a mean-field approximation \cite{biswal2020stabilization} of the state evolution of the agents. The problem of steering the group of agents is then turned to steering an initial density function to a terminal one. The density control problem is formulated as follows.

\begin{problem}
Consider the system dynamics
\begin{equation}
x(k+1) = a(k)x(k) + u(k), 
\label{systemeq}
\end{equation}
where the random variables $x(k)$ and $u(k)$ take values in $\mathbb{R}$, and $\{a(k)\}_{k=0}^{K-1}$ is a deterministic sequence.

The group steering problem is formulated as given an initial probability density function $q_{0}(x) \in \mathcal{D}_{2n}$ for the random variable $x(0)$, determine a control sequence $\left( u(0), \cdots, u(K-1) \right)$ such that the terminal state $x(K)$ has the prescribed density $\tau(x) \in \mathcal{D}_{2n}$.
\end{problem}

However it is not always feasible for us to obtain a closed form of solution to this problem. If the distributions are not assumed to belong to parametric families such as the Gaussian or exponential family, the problem is intrinsically infinite-dimensional. Denote $f_{x(k), u(k)}\left(t , \xi \right)$ as joint distribution of $x(k)$ and $u(k)$. We note that the density function of $x(k+1)$ can be written as
\begin{equation}
\begin{aligned}
    q_{k+1}(t) & = \int_{\mathbb{R}} f_{x(k), u(k)}\left( t, \xi \right)d\xi\\
    & = \frac{1}{a(k)}\int_{\mathbb{R}} q_{k}\left(\frac{\xi}{a(k)}\right) p_{k}\left(t-\xi\right) d\xi\\
    & = \frac{1}{a(k)}\left(q_k\left(\frac{s}{a(k)}\right) * p_k(s)\right)(t).
\end{aligned}
\label{qk1}
\end{equation}
where in the second equality we have used the fact that $x(k)$ and $u(k)$ are independent. Addressing the density steering problem requires obtaining an analytical solution for $q_{k+1}(t)$ in \eqref{qk1}. The scenario would be more complicated with $x(k), u(k)$ being correlated. However, except for specific function classes like Gaussian distributions or trigonometric functions, obtaining such an analytic solution is impractical for the majority of function classes. This limitation is a key factor, as evident in prior results with a problem setting akin to the density steering problem, where the examples predominantly feature Gaussian or trigonometric densities. Consequently, the applicability of these results in real-world scenarios is severely limited.

There is a similar problem in non-Gaussian Bayesian filtering. In our previous results \cite{wu2023non}, we proposed a method of using the power moments to treat this intractable problem, mainly for characterizing the macroscopic property of the distributions. However, even it is theoretically feasible to characterize the distribution of the agents by the full power moment sequence, the problem is infinite dimensional. A common treatment is to truncate the first $2n$ moment terms \cite{byrnes2003convex, georgiou2003kullback}, which turns the problem we treat to a truncated moment problem.

By the system equation \eqref{systemeq}, the power moments of the states up to order $2n$ are written as $
\mathbf{E}\left[ x^{l}(k+1) \right] = \sum_{j=0}^{l}\binom{l}{j}a^{j}(k)\mathbf{E}\left[ x^{j}(k)u^{l-j}(k) \right]$ for $l \in \mathbb{N}_{0}$ ($\mathbb{N}_{0}$ denotes the set of all nonnegative integers), $l \leq 2n$. It is difficult to treat the term $\mathbf{E}\left[ x^{j}(k)u^{l-j}(k) \right]$. However, we note that if $x(k)$ and $u(k)$ are independent from each other, i.e., $\mathbf{E}\left[ x^{j}(k)u^{l-j}(k) \right] = \mathbf{E}\left[ x^{j}(k) \right]\mathbf{E}\left[ u^{i-j}(k) \right]$, the dynamics of the moments can be written in a linear matrix equation
\begin{equation}
    \mathscr{X}(k+1) = \mathscr{A}(\mathscr{U}(k))\mathscr{X}(k)+\mathscr{U}(k)
\label{momentsystem}
\end{equation}
where the new state vector is composed of the power moment terms up to order $2n$,
\begin{equation}
\mathscr{X}(k) = \begin{bmatrix}
\mathbf{E}[x(k)] & \mathbf{E}[x^{2}(k)] & \cdots & \mathbf{E}[x^{2n}(k)]
\end{bmatrix}^{\intercal}.
\label{XK}
\end{equation}

The new input vector is written as
\begin{equation}
\mathscr{U}(k) = \begin{bmatrix}
\mathbf{E}[u(k)] & \mathbf{E}[u^{2}(k)] & \cdots & \mathbf{E}[u^{2n}(k)]
\end{bmatrix}^{\intercal}.
\label{UK}
\end{equation}
where we have denoted
\begin{equation}
\mathbf{E}\left[ x^{l}(k) \right] = \int_{\mathbb{R}}x^{l}q_{k}(x)dx,
\label{xlK}
\end{equation}
and
\begin{equation}
\mathbf{E}\left[ u^{l}(k) \right] = \int_{\mathbb{R}}u^{l}p_{k}(u)du.
\label{ulK}
\end{equation}
The matrix $\mathscr{A}(\mathscr{U}(k))$ of the new system can then be written as \eqref{longeq1}. By employing truncated power moments to describe the dynamics of the system \eqref{uniequation}, where $x(k)$ and $u(k)$ represent probability densities, we aim to reformulate the control problem as the manipulation of the power moments of these probability densities. The system \eqref{momentsystem} is denoted as the moment system associated with \eqref{uniequation}.

\begin{figure*}[t]
\begin{equation}
\mathscr{A}(\mathscr{U}(k))
= \begin{bmatrix}
a(k) & 0 & 0 & \cdots & 0\\ 
2a(k)\mathbf{E}[u(k)] & a^{2}(k) & 0 & \cdots & 0\\ 
3a(k)\mathbf{E}[u^{2}(k)] & 3a^{2}(k)\mathbf{E}[u(k)] & a^{3}(k) & \cdots & 0\\ 
\vdots & \vdots & \vdots & \ddots\\ 
\binom{2n}{1}a(k)\mathbf{E}[u^{2n-1}(k)] & \binom{2n}{2}a^{2}(k)\mathbf{E}[u^{2n-2}(k)] & \binom{2n}{3}a^{3}(k)\mathbf{E}[u^{2n-3}(k)] &  & a^{2n}(k)
\end{bmatrix}.
\label{longeq1}
\end{equation}
\hrulefill
\vspace*{4pt}
\end{figure*}

By the proposed moment system, the original density steering problem can be reduced to distribution steering by power moments, which is formulated as follows.
\begin{problem}
The dynamics of the moment system is \eqref{momentsystem}, where $\mathscr{X}(k), \mathscr{U}(k)$ are obtained by \eqref{xlK},\eqref{ulK}. Given an arbitrary initial density $q_{0}(x) \in \mathcal{D}_{2n}$, determine the control sequence $\left( u(0), \cdots, u(K-1)\right)$ such that the first $2n$ order power moments of the terminal density are identical to those of a specified one, i.e.,
\begin{equation}
\int_{\mathbb{R}} x^{l} \tau(x) dx = \int_{\mathbb{R}} x^{l} q_{K}(x) dx
\label{ExTl}
\end{equation}
for $l = 1, \cdots, 2n$, where $\tau \in \mathcal{D}_{2n}$ is the specified terminal density function.
\label{Def22}
\end{problem}
\begin{remark}
The advantage of the proposed problem formulation is twofold. Firstly, with this formulation, the problem becomes finite-dimensional, allowing for a closed-form solution. Secondly, it eliminates the need for the initial and terminal density functions of the agents to adhere to specific parametric families, enhancing the algorithm's applicability across a diverse array of real-world scenarios.
\end{remark}


The task is now to figure out an algorithm to determine a sequence of $\left( \mathscr{U}(0), \cdots, \mathscr{U}(K-1) \right)$. Nonetheless, there exist two primary distinctions from conventional control problems. Firstly, in this problem, the system matrix of the moment system becomes a function of the control vector. Secondly, the sequence of elements in the control vector, denoted as $\mathscr{U}(k)$, must adhere to the condition that the corresponding Hankel matrix, namely
$$
\begin{aligned}
\mathbf{H}_{\mathscr{U}(k)}=\mathbf{E}\left[\left[\begin{array}{c}u^0(k) \\ u^1(k) \\ \vdots \\ u^n(k)\end{array}\right]\left[\begin{array}{llll}u^0(k) & u^1(k) & \cdots & u^n(k)\end{array}\right]\right]
\end{aligned}
$$
is positive definite. It is a classic result of the truncated Hamburger moment problem (see e.g. \cite[Corollary 9.2]{schmudgen2017moment}), which concerns the characterization of measures defined on $\mathbb{R}$ that satisfy a given sequence of moment conditions. Here $\boldsymbol{\mathrm{H}}_{\mathscr{U}}$ denotes the Hankel matrix of the vector $\mathscr{U}$. We define such a subspace of $\mathbb{R}^{2n}$ as $\mathbb{V}^{2n}_{++} := \{ \mathscr{U} \in \mathbb{R}^{2n} \mid \boldsymbol{\mathrm{H}}_{\mathscr{U}(k)} \succ 0 \}$.

In prior studies, the distribution steering problems have consistently employed optimal control strategies. However, such an approach is not entirely feasible in this specific problem. The challenge lies in the necessity to ensure that both $\mathscr{X}(k)$ and $\mathscr{U}(k)$ belong to $\mathbb{V}^{2n}_{++}$. To the best of our knowledge, there has yet to be a viable result capable of addressing the optimal control problem while imposing the constraint that states and control inputs must reside within $\mathbb{V}^{2n}_{++}$ - in other words, ensuring that the corresponding Hankel matrices are positive definite.

Let $\mathcal{U}$ be the feasible set of control sequences $\boldsymbol{\mathscr{U}} := \left(\mathscr{U}(0), \cdots, \mathscr{U}(K-1) \right)$, which satisfies $
\sum_{k=0}^{K-1}\mathscr{U}^{\intercal}(k)\mathscr{U}(k) < \infty$ and effects the terminal system state $x(K)$ to be distributed satisfying \eqref{ExTl}. Then the family $\mathcal{U}$ represents admissible control inputs which achieve the desired moment transfer.
Denote the error of the moments from the specified ones as
\begin{equation}
    e(k) = \mathscr{X}_{T} - \mathscr{X}(k),
\label{ek}
\end{equation}
where the elements of $\mathscr{X}_{T}$ are the power moments corresponding to the specified terminal density function $\tau(x)$.

Simultaneously confining both $\mathscr{U}(k)$ and $\mathscr{X}(k)$ to belong to the set $\mathbb{V}^{2n}_{++}$ is not always feasible. However, it's worth noting that a sub-optimal solution to Problem \ref{Def22} can be attained by first determining the trajectory of the state and subsequently obtaining the control inputs corresponding to this trajectory. A similar approach was utilized in \cite{sivaramakrishnan2022distribution} to address the general distribution steering problem. We first determine the state trajectory. Define the initial error as $e(k_{0}) := \mathscr{X}_{T} - \mathscr{X}(k_{0})$, and express the state at step $k$ as
\begin{equation}
\mathscr{X}(k) = \mathscr{X}(k_{0}) + \left( \sum_{i = k_{0}}^{k - 1} \omega_{i} \right) e(k_{0}),
\label{Xksum}
\end{equation}
where the weighting coefficients $\{\omega_{i}\}$ to be determined satisfy $\omega_{i} \in \mathbb{R}_{+}$ for $i = k_{0}, \dots, K-1$ and $\sum_{i = k_{0}}^{K - 1} \omega_{i} = 1$. Then we have the following lemma.
\medskip

\begin{lemma}
Given $e(k_{0})\in \mathbb{V}^{2n}_{++}$, we have $\mathscr{X}(k) \in \mathbb{V}^{2n}_{++}$
for $k = k_{0}+1, \cdots, K$ where $\omega_{i} \in \mathbb{R}_{+}$ for $i = k_{0}, \cdots, K-1$ and $\sum_{i = k_{0}}^{K - 1} \omega_{i} = 1$. 
\label{Lemma21}
\end{lemma}
\begin{proof}
The proof is straightforward. Since $\mathscr{X}(k_{0}), e(k_{0}) \in \mathbb{V}^{2n}_{++}$, we have $\boldsymbol{\mathrm{H}}_{\mathscr{U}(k)} \succ 0$, $\boldsymbol{\mathrm{H}}_{e(k_{0})} \succ 0$. We note that the sum of positive definite matrices is still positive definite. Since $\omega_{i} > 0$, we have $\omega_{i} e(k_{0}) \in \mathbb{V}^{2n}_{++}$. Then $\mathscr{X}(k) \in \mathbb{V}^{2n}_{++}$.
\end{proof}

Now it remains to prove that there exists a time step $k_{0}$ at which $\mathscr{X}_{T} - \mathscr{X}(k_{0}) \in \mathbb{V}^{2n}_{++}$.
\medskip

\begin{proposition}
There exists a time step $k_{0}$ which satisfies $\mathscr{X}_{T} - \mathscr{X}(k_{0}) \in \mathbb{V}^{2n}_{++}$, assuming that $\mathscr{X}(k)$, $0 \leq k \leq k_{0}$, are uncontrolled states, i.e., $u(k) = 0$ for $0 \leq k \leq k_{0}$, and that $\alpha_k := \prod_{i=0}^{k-1} a(i) \to 0$ as $k\to\infty$ (equivalently, $\sum_{k=0}^{\infty} (1 - a(k)) = \infty$).
\label{proposition22}
\end{proposition}

\begin{proof}
We need to prove that there exists $k_0 \in \mathbb{N}_0$ such that $\boldsymbol{\mathrm{H}}_{\mathscr{X}_T - \mathscr{X}(k_0)} \succ 0$, with $u(k) = 0$ for $0 \leq k \leq k_0$. By the definition of positive definiteness, this is equivalent to proving that each leading principal minor, i.e., the determinant of each leading principal submatrix, is positive.

Let $\alpha_k := \prod_{i=0}^{k-1} a(i) \in (0,1)$. Denote the $i$-th order leading principal submatrix of $\boldsymbol{\mathrm{H}}_{\mathscr{X}_T - \mathscr{X}(k)}$ as $H_i(k)$. Define the real-valued function $f_i : (0,1) \to \mathbb{R}$, $f_i(\alpha_k) := \det(H_i(k))$, $i = 1,\dots,2n$. Moreover, we note that each $f_i(\cdot)$ is a real polynomial with even degree. Since the system is uncontrolled and $\alpha_k \to 0$, it follows from \eqref{momentsystem} and \eqref{longeq1} that $\mathscr{X}(k) \to 0$ as $k \to \infty$. Hence $\boldsymbol{\mathrm{H}}_{\mathscr{X}_T - \mathscr{X}(k)} \to \boldsymbol{\mathrm{H}}_{\mathscr{X}_T} \succ 0$. Therefore, each polynomial $f_i(\alpha_k) \to f_i(0) > 0$ as $k \to \infty$.

Since each $f_i(\cdot)$ is continuous on $[0,1]$, for $i = 1,\dots,2n$, there exists $\epsilon_i > 0$ such that $f_i(\alpha_k) = \det(H_i(k)) > 0$ whenever $0 < \alpha_k < \epsilon_i$. Now, since $\alpha_k \to 0$, there exists an integer $k_i$ such that $\alpha_k < \epsilon_i$ for all $k \geq k_i$. Define $k_0 := \max_i k_i$. Then for all $i = 1,\dots,2n$, we have $\det(H_i(k_0)) > 0$, which implies that all leading principal minors of $\boldsymbol{\mathrm{H}}_{\mathscr{X}_T - \mathscr{X}(k_0)}$ are strictly positive. Hence, $\boldsymbol{\mathrm{H}}_{\mathscr{X}_T - \mathscr{X}(k_0)} \succ 0$, i.e., $\mathscr{X}_T - \mathscr{X}(k_0) \in \mathbb{V}^{2n}_{++}$, completing the proof.
\end{proof}

\begin{remark}
The additional assumption that $\alpha_k := \prod_{i=0}^{k-1} a(i) \to 0$ as $k\to\infty$
is mild and satisfied by a broad class of sequences arising in dynamical systems. For instance, if $a(k)$ takes values in any closed subset of $(0,1)$ that is strictly bounded away from $1$, then $\alpha_k \to 0$ trivially. In practical swarm applications (e.g., UAVs, microsatellites), the discrete-time dynamics $x_i(k+1) = a(k)x_i(k) + u_i(k)$ typically result from the zero-order hold discretization of a stable continuous-time system $\dot{x}_i(t) = -\beta(t)x_i(t) + u_i(t)$. The parameter $\beta(t) > 0$ represents inherent physical dissipation (e.g., aerodynamic drag) or a low-level stabilizing controller, ensuring the uncontrolled open-loop dynamics ($u_i(t) = 0$) are asymptotically stable \cite{beard2012small}. With a sampling period $\Delta t > 0$, the discrete system parameter is exactly $a(k) = e^{-\beta(k)\Delta t}$. For such physical systems, it is a common property that the damping coefficient is bounded away from zero, i.e., $\beta(k) \ge \beta_{\min} > 0$. It immediately follows that $a(k) = e^{-\beta(k)\Delta t} \le e^{-\beta_{\min}\Delta t} := 1 - \epsilon < 1$. Because $a(k)$ is strictly bounded away from 1 by a fixed margin $\epsilon > 0$, the sequence $\alpha_k = \prod_{i=0}^{k-1} a(i)$ satisfies $\alpha_k \le (1-\epsilon)^k$, converging to 0 exponentially fast. This physically grounded bound rigorously justifies the assumption in Proposition~\ref{proposition22}.
\end{remark}

By Proposition \ref{proposition22}, it is possible to choose a time step $k_{0}$ which satisfies $\mathscr{X}_{T} - \mathscr{X}(k_{0}) \in \mathbb{V}^{2n}_{++}$. We assume that the system is uncontrolled before $k_{0}$, i.e. $u(k)=0, k \leq k_{0}$. From step $k_{0}$, we impose controls on the system. Lemma \ref{Lemma21} has proved the positiveness of $\mathscr{X}(k), k = k_{0}, \cdots, K$. Therefore it remains to determine the parameters $\omega_{k}, k = k_{0}, \cdots, K - 1$ and the corresponding control inputs $\mathscr{U}(k)$.

It is a non-trivial problem. We give an empirical scheme to treat it. To ensure gradual evolution of the state trajectory in discrete time, it is desired that the $\omega_{i}$'s are close to each other. It is usually feasible for us to choose $
\omega_{k_{0}} = \cdots = \omega_{K-1} = \frac{1}{K-k_{0}}$. We note that to determine the trajectory of the system state \textit{a priori} has also been adopted in the literature. In \cite{sivaramakrishnan2022distribution}, prespecified system states, which was called to form a reference trajectory, were used to determine the control inputs. After that the parameters $\omega_{i}'s$ are determined, the control inputs of the moment system $\mathscr{U}(i)$ for $i = k_{0}, \cdots, K-1$ can then be calculated by solving the equation \eqref{momentsystem}, provided with $\mathscr{X}(k), k = k_{0}+1, \cdots, K$ calculated by \eqref{Xksum}.

However, in some cases, the control inputs \( \mathscr{U}(k) \) may fail to lie in \( \mathbb{V}^{2n}_{++} \) when the weights \( \omega_i \) are chosen to be equal. This typically occurs when the specified terminal density is multi-modal (i.e., has multiple peaks). For such densities, the absolute values of higher-order even moments — such as \( \left|\mathbf{E}[x^{2n}(K-1)]\right| \) and \( \left|\mathbf{E}[x^{2n-2}(K-1)]\right| \) — tend to be small. As a result, the Hankel matrix \( H_{\mathscr{U}(K-1)} \) is more likely to not be positive definite (\( H_{\mathscr{U}(K-1)} \nsucc 0 \)), leading to non-existence of $u(K-1)$. If so, we can choose a larger $\omega_{K-1}$ and let 
\begin{equation}
\omega_{k_{0}} = \cdots = \omega_{K-2} = \frac{1 - \omega_{K-1}}{K-k_{0}}.
\label{Adjustk}
\end{equation}

It should be noted that this approach is empirical, based on numerical evidence and not on a formal proof. Other selections of $\omega_k$ may be employed to prescribe the trajectory \textit{a priori}. Now we need to obtain the $u(k)$ given the power moments of the controls for the moment system, which we call realization of the control inputs. This problem falls within the general framework of the Hamburger moment problem \cite{schmudgen2017moment}.  Inspired by \cite{georgiou2003kullback}, we propose to treat this problem by minimizing the Kullback-Leibler distance between a reference density $r(u)$ and the probability density $p(u)$, with the moment constraints. For the sake of simplicity, we omit $k$ if there is no ambiguity in the following part of this section.

The Kullback--Leibler distance is defined as
\begin{equation}
    \mathrm{KL}(r \,\|\, p)
    = \int_{\mathbb{R}} r(u)\,\log\frac{r(u)}{p(u)} \, \mathrm{d}u,
    \label{KL-def}
\end{equation}
where $r$ is an arbitrary reference probability density in $\mathcal{D}$.
Let $\mathbb{S}^{n+1}:= \{ X \in \mathbb{R}^{(n+1)\times(n+1)} \mid X^\top = X \}$ be the space of real symmetric $(n+1)\times(n+1)$ matrices, and let $B(t) :=
\begin{bmatrix}
1 & t & \cdots & t^{n-1} & t^{n}
\end{bmatrix}^{\!\top}$.
We introduce the linear space of signed measures
\[
    \mathcal{F}
    := \left\{ \mu : \int_{\mathbb{R}} |t|^k \mathrm{d}|\mu|(t) < \infty, \quad k=0, \ldots, 2n \right\},
\]
and define the linear integral operator $\Theta: \mathcal{F} \rightarrow \mathbb{S}^{n+1}$ by
\begin{equation}
    \Theta(\mu)
    := \int_{\mathbb{R}} B(t)B(t)^{\top} \,\mathrm{d}\mu(t).
\end{equation}
The set $\mathcal{D}_{2n}$ denotes the set of probability density functions with finite moments up to order $2n$. We identify $\mathcal{D}_{2n}$ with a convex subset of $\mathcal{F}$ via the correspondence $\mathrm{d}\mu(t) = p(t)\,\mathrm{d}t$. In the rest of the paper, we abuse notation slightly to write $\Theta(p)$ for the integral restricted to this density, i.e., $\Theta(p) := \int B(t)B(t)^\top p(t)\mathrm{d}t$. We also note that
$\operatorname{range}(\Theta) := \Theta(\mathcal{F})$ is a linear subspace of
$\mathbb{S}^{n+1}$, while $\Theta(\mathcal{D}_{2n})$ is a convex subset of
$\operatorname{range}(\Theta)$.

\begin{theorem}[Theorem~3.2, \cite{wu2023non}]
Let $\mathcal{L}_{+} := \left\{\Lambda \in \operatorname{range}(\Theta)\ \Big|\ B(u)^{\top} \Lambda B(u) > 0,\ \forall u \in \mathbb{R}\right\}$. Since $\mathcal{D}_{2n}$ is convex and $\Theta$ is linear on the embedding space $\mathcal{F}$,
$\Theta(\mathcal{D}_{2n})$ is a convex subset of $\operatorname{range}(\Theta)$.
For any $r \in \mathcal{D}$ and any positive definite Hankel matrix
$\Sigma \in \Theta(\mathcal{D}_{2n})$ with $\Sigma(1,1) = 1$, there exists a unique
$\hat{p} \in \mathcal{D}_{2n}$ that minimizes \eqref{KL-def} subject to
$\Theta(\hat{p}) = \Sigma$. Specifically,
\begin{equation}
    \hat{p}(u)
    = \frac{r(u)}{B(u)^{\top} \hat{\Lambda} B(u)},
    \label{hatp}
\end{equation}
where $\hat{\Lambda}$ is the unique solution to the minimization problem
\begin{equation}
    \mathbb{J}_{r}(\Lambda)
    :=
    \operatorname{tr}(\Lambda \Sigma)
    - \int_{\mathbb{R}} r(u)\,\log \big( B(u)^{\top} \Lambda B(u) \big)\,
        \mathrm{d}u
    \label{Jr}
\end{equation}
over all $\Lambda \in \mathcal{L}_{+}$.
Here $\operatorname{tr}(\cdot)$ denotes the trace of a matrix.
\end{theorem}
Then the density estimation is formulated as a convex optimization problem. Unlike other methods of moments, the power moments of our proposed density estimate are exactly identical to those specified, which makes it a satisfactory approach for realization of the control inputs. Since the prior density $r(u)$ and the density estimate $\hat{p}(u)$ are both supported on $\mathbb{R}$, $r(u)$ can be chosen as a Gaussian distribution (or a Cauchy distribution if $\hat{p}(u)$ is assumed to be heavy-tailed).

We conclude the algorithm for density steering of the vast group of agents in this section as in Algorithm \ref{alg:1}.

\medskip

\begin{algorithm}
    \caption{Density Steering of a Group of Agents.}
    \label{alg:1}
    \begin{algorithmic}[1]
        \Require The maximal time step $K$; the parameter of the system $a(k)$ for $k = 0, \cdots, K-1$; the initial system density $q_{0}(x)$; the specified terminal density $\tau(x)$.
        \Ensure The controls $u(k)$, $k = 0, \cdots, K-1$.
        \State $k \Leftarrow 0$
    \While{$k < K$ and $e(k) \notin \mathbb{V}^{2n}_{++}$}
        \State Calculate $\mathscr{X}(k)$ by \eqref{momentsystem} if $k > 0$ or by \eqref{XK} if $k = 0$
        \State Calculate $e(k)$ by \eqref{ek}
        \If{$e(k) \in \mathbb{V}^{2n}_{++}$}
        \State Calculate the states of the moment system $\mathscr{X}(i)$ for $i = k + 1, \cdots, K-1$ by \eqref{Xksum} with $\omega_{k} = \cdots = \omega_{K-1}$ \label{Step6}
        \State Calculate the controls of the moment system $\mathscr{U}(i)$ for $i = k, \cdots, K-1$ by \eqref{momentsystem}
        \If{$\exists i, \mathscr{U}(i) \notin \mathbb{V}^{2n}_{++}$}
        \State Back to Step \ref{Step6}, increase $\omega_{K-1}$ and determine other $\omega_{0}, \cdots, \omega_{K-2}$ by \eqref{Adjustk}
        \EndIf
        \State Optimize the cost function \eqref{Jr} and obtain the analytic estimates of the densities $\hat{p}_{i}(u)$ for $i = k, \cdots, K-1$
        \Else
            \State $u(k) = 0$
        \EndIf
        \State Calculate the power moments of the system state $x(k+1)$, i.e., $\mathscr{X}(k+1)$
        \State $k \Leftarrow k+1$
    \EndWhile
    \end{algorithmic}
\end{algorithm}

Through the proposed density steering scheme, the first $2n$-order power moments of the terminal density function match those of the desired one. As $n$ approaches infinity, the full moment sequence of the terminal density converges to that of the desired density. By Theorem 4.5.5 in \cite{chung2001course}, we can straightforwardly demonstrate that the terminal distribution obtained from our distribution steering algorithm converges almost everywhere to the true desired distribution as $n$ tends to infinity. Assuming both the initial and desired terminal distributions are continuous, we can further conclude that the terminal density generated by the proposed algorithm converges to the desired distribution. This convergence to the desired density distinguishes our proposed algorithm from other existing methods. Since we used the truncated power moments of the initial and terminal density functions for steering, there may exist an error between the terminal density and the desired one. We will propose an upper bound of error of the terminal density, as to characterize the maximal difference between the terminal density by our proposed algorithm and the specified one. In \cite{wu2023non}, we proposed an error upper bound for the Hamburger moment problem in the sense of Kolmogorov-Smirnov distance, which is a measure widely used in the moment problem \cite{Aldo2003A, tagliani2003maximum}.

By our proposed density steering scheme, $\tau$ and $q_K$ have the same set of power moments up to order $2n$. The Kolmogorov-Smirnov distance between the terminal density $q_{K}(x)$ and the desired terminal density $\tau(x)$ is defined as $V(q_{K}, \tau) = \sup_{t} \left|\int_{\left(-\infty, t \right] }(q_{K} - \tau) d t\right| = \sup_{t} \left| F_{q_{K}} - F_{\tau} \right|$, where $F_{q_{K}}$ and $F_{\tau}$ are the two cumulative distribution functions of $q_{K}$ and $\tau$.

Differential-entropy is used to calculate the upper bound of Kolmogorov-Smirnov distance in \cite{Aldo2003A}. The differential-entropy \cite{michalowicz2013handbook} is defined as $H[q] = - \int_{\mathbb{R}}q(t) \log q(t)dt$. We first introduce the differential-entropy maximizing distribution $F_{\breve{q}_{K,2n}}$, of which the moments are the power moments of $q_{K}$. It has the density function \cite{kapur1992entropy} $
    \breve{q}_{K,2n}(t) = \exp \left ( - \sum_{i = 0}^{2n} \lambda_{i} t^{i} \right ) $ where $\lambda_{0}, \cdots, \lambda_{2n}$ are determined by the following constraints, namely $
    \int_{\mathbb{R}} t^k \exp \left ( - \sum_{i = 0}^{2n} \lambda_{i} t^{i} \right )dt=\int_{\mathbb{R}} t^k \tau(t) dt$, for $k=0,1, \cdots, 2n$. 
\begin{proposition}
    The Kolmogorov-Smirnov distance between the desired terminal density function $\tau(x)$ and the terminal density by our proposed steering scheme $q_{K}(x)$, is upper bounded by
\begin{equation}
\begin{aligned}
& V \left ( q_{K}, \tau \right ) \\
\leq & 3 \left[-1+\left\{1+\frac{4}{9}\left(H\left [ \breve{q}_{K,2n} \right ] - H\left [ q_{K} \right ]\right)\right\}^{1 / 2}\right]^{1 / 2} \\
+ & 3\left[-1+\left\{1+\frac{4}{9}\left(H\left [ \breve{q}_{K,2n} \right ] - H\left [ \tau \right ]\right)\right\}^{1 / 2}\right]^{1 / 2}.
\end{aligned}
\label{Vqktau}
\end{equation}
\end{proposition}

\begin{proof}
We also note that 
$$
\begin{aligned}
H[\breve{q}_{K,2n}] = & - \int_{\mathbb{R}}\breve{q}_{K,2n}(t) \log \breve{q}_{K,2n}(t)dt\\
= & \sum_{i = 0}^{2n} \lambda_{i} \int_{\mathbb{R}}t^{i}\breve{q}_{K,2n}(t)dt.
\end{aligned}
$$
By referring to \cite{Aldo2003A}, the KL distance between the true density and the differential-entropy maximizing density can be written as
$$
\begin{aligned}
    KL \left(q_{K}\| \breve{q}_{K,2n}\right) = & \int_{\mathbb{R}} q_{K}(t) \log \frac{q_{K}(t)}{\breve{q}_{K,2n}(t)} d t\\
    = & - H\left [ q_{K} \right ] + \sum_{i = 0}^{2n} \lambda_{i} \int_{\mathbb{R}} t^i \breve{q}_{K,2n}(t) dt\\
    = & H\left [ \breve{q}_{K,2n} \right ] - H\left [ q_{K} \right ].
\end{aligned}
$$
Similarly, we can obtain 
$$
\begin{aligned}
KL \left(\tau\| \breve{q}_{K,2n}\right) = & \int_{\mathbb{R}} \tau(t) \log \frac{\tau(t)}{\breve{q}_{K,2n}(t)} d t\\
    = & - H\left [ \tau \right ] + \sum_{i = 0}^{2n} \lambda_{i} \int_{\mathbb{R}} t^i \tau(t) dt\\
    = & H\left [ \breve{q}_{K,2n} \right ] - H\left [ \tau \right ].
\end{aligned}
$$

By \cite{1970Correction, Aldo2003A}, we obtain
\begin{equation}
\begin{aligned}
    & V \left ( \breve{q}_{K,2n}, q_{K} \right )\\
    \leq & 3\left[-1+\left\{1+\frac{4}{9} KL \left(q_{K} \| \breve{q}_{K,2n} \right)\right\}^{1 / 2}\right]^{1 / 2} \\
    = & 3\left[-1+\left\{1+\frac{4}{9} \left ( H\left [ \breve{q}_{K,2n} \right ] - H\left [ q_{K} \right ] \right )\right\}^{1 / 2}\right]^{1 / 2}
\end{aligned}
\label{Vbound1}
\end{equation}
and
\begin{equation}
\begin{aligned}
    & V \left ( \breve{q}_{K,2n}, \tau \right )\\
    \leq & 3\left[-1+\left\{1+\frac{4}{9} \left ( H\left [ \breve{q}_{K,2n} \right ] - H\left [ \tau \right ] \right )\right\}^{1 / 2}\right]^{1 / 2}.
\end{aligned}
\label{Vbound2}
\end{equation}
Then the error upper bound of the terminal density can be written as
\begin{equation}
\begin{aligned}
& V \left ( q_{K}, \tau \right ) \\
= & \sup_{t}|F_{q_{K}}\left ( t \right )-F_{\tau}\left ( t \right )| \\
\leq & \sup_{t} \left|F_{q_{K}}\left ( t \right )-F_{\breve{q}_{K,2n}}\left ( t \right )\right|+\sup_{t} \left|F_{\breve{q}_{K,2n}}(t)-F_{\tau}(t)\right|\\
= & V \left ( \breve{q}_{K,2n}, q_{K} \right ) + V \left ( \breve{q}_{K,2n}, \tau \right ).
\end{aligned}
\label{UpperBoundUnbiased}
\end{equation}
By substituting \eqref{Vbound1} and \eqref{Vbound2} into \eqref{UpperBoundUnbiased}, we obtain \eqref{Vqktau}, which completes the proof.
\end{proof}

\begin{remark}
We note that $\breve{q}_{K,2n}$ is obtained by maximizing the differential entropy subject to moment constraints up to order $2n$, and the bound $V(q_K,\tau)$ depends on $2n$. Let $\mathcal{M}_{2n} \subset \mathcal{D}_{2n}$ denote the set of all probability density functions satisfying the moment constraints induced by $\tau$ up to order $2n$. It is clear that $\mathcal{M}_{2n_2} \subset \mathcal{M}_{2n_1}$ whenever $n_2 > n_1$, since imposing more moment constraints reduces the feasible set.
Consequently, by monotonicity of maximization over nested feasible sets, the entropy of the maximum-entropy solution,
$H[\breve{q}_{K,2n}] = \max_{p\in\mathcal{M}_{2n}} H(p)$,
is non-increasing as $2n$ increases.
Moreover, since both $q_K$ and $\tau$ belong to $\mathcal{M}_{2n}$, we have
$H[\breve{q}_{K,2n}] \ge H[q_K]$ and $H[\breve{q}_{K,2n}] \ge H[\tau]$.
Therefore, by \eqref{Vqktau}, the resulting upper bound on $V(q_K,\tau)$ is non-increasing in $2n$ (i.e., it becomes tighter as $2n$ increases).
We also observe that the first $2n$ moments of $\breve{q}_{K,2n}$, $q_K$, and $\tau$ coincide.
Under standard moment determinacy conditions \cite{chung2001course}, as $2n\to\infty$, the distributions
$\breve{q}_{K,2n}$, $q_K$, and $\tau$ coincide almost everywhere.
As a consequence, both
$H[\breve{q}_{K,2n}] - H[q_K]$ and
$H[\breve{q}_{K,2n}] - H[\tau]$
tend to zero, which implies that $V(q_K,\tau)\to 0$ as $n\to\infty$.
\end{remark}

\section{Steer the group as an occupation measure}

In the preceding section, we introduced an algorithm designed to guide a vast group of agents, represented as a probability density function, to a specified terminal state. However, the characterization of agents as a density function is an approximation when the number of agents tends to infinity. To apply our proposed algorithm to real-world scenarios, where the number of agents is finite, we must address individual agents and provide specific control inputs for each. In this section, we employ the occupation measure to describe the agent group and present a control scheme aiming at steering an initial occupation measure towards a terminal one. We first define the occupation measure \cite{zhang2020modeling} of the agents at time step $k$ by 
$$
d \mu_{k}(x, u)=\frac{1}{N} \sum_{i=1}^{N} \delta\left(x-x_{i}(k)\right) \delta\left(u-u_{i}(k)\right) d x d u,
$$
where $N$ is the number of the agents. Then the occupation measure of the group state can be written as $ dq_{k}(x) = \frac{1}{N}\sum_{i=1}^{N} \delta\left(x-x_{i}(k)\right)dx$, 
which is a marginal measure of $d \mu_{k}(x, u)$. Another marginal measure, the occupation measure of the controls on the group of agents reads $dp_{k}(u) = \frac{1}{N}\sum_{i=1}^{N} \delta\left(u-u_{i}(k)\right)du$.

Now we define the occupation measure steering by power moments, which we will treat in this section.

\begin{problem}
The dynamics of the moment system is \eqref{momentsystem}, where $\mathscr{X}(k), \mathscr{U}(k)$ are defined as $\mathbf{E}\left[x^{l}(k) \right] = \int_{\mathbb{R}}x^{l}dq_{k}(x) = \frac{1}{N}\sum_{i=1}^{N}x^{l}_{i}(k)$ and $\mathbf{E}\left[ u^{l}(k) \right] = \int_{\mathbb{R}}u^{l}dp_{k}(u) = \frac{1}{N}\sum_{i=1}^{N}u^{l}_{i}(k)$. Given an arbitrary initial state $x(0)$, determine the control sequence $\left( u_{i}(0), \cdots, u_{i}(K-1)\right)$ for each agent $i$ such that the power moments of the terminal occupation measure are identical to those as specified, i.e., $\mathbf{E}\left[ x_{T}^{l} \right] = \int_{\mathbb{R}} x^{l} \tau(x) dx = \frac{1}{N} \sum_{i=1}^{N} x_{i}^{l}(K)$ for $l = 1, \cdots, 2n$, where $\tau(x)$ is the terminal density.
\end{problem}

The main difference of the occupation measure steering problem from the density steering one lies in determining the control input for each agent. We naturally consider designing feedback control laws for the agents. It might be feasible with a limited number of agents. However it is quite expensive and problematic with quite a large number of agents as in our problem. Since $N$ is large, we consider first estimating the occupation measure of $u(k)$ as a continuous function $\hat{p}_{k}(u)$. Then we draw $N$ i.i.d samples from it and assign them to $u_{i}(k)$, i.e., $u_{i}(k) \sim \hat{p}_{k}(u), i \in \mathbb{N}_{0}, i \leq N$. By the strong law of large numbers, we note that $\frac{1}{N}\sum_{i=1}^{N}u^{l}_{i}(k) \stackrel{a.s.}{\longrightarrow} \int_{\mathbb{R}} u^{l}\hat{p}_{k}(u)du, \text{with} \ N \rightarrow +\infty$, 
which means that the power moments of $u(k)$ converge almost surely to the power moments $\mathscr{U}(k)$ of the designed controls. Moreover, the sampling strategy ensures that the system state $x(k)$ and the control input $u(k)$ are independent from each other. Then the problem comes to putting forward a sampling strategy. We consider using the acceptance-rejection sampling \cite{brandimarte2014handbook} strategy. 

We now implement the proposed acceptance-rejection sampling strategy to modify Algorithm \ref{alg:1} for addressing the occupation measure steering problem, resulting in Algorithm \ref{alg:3}. To the best of our knowledge, Algorithm \ref{alg:3} is the first scalable general distribution steering algorithm capable of handling a large number of agents, including scenarios with thousands of agents or more. Moreover, we could approximate $H\left [ q_{K} \right ]$ by $ H\left [ q_{K} \right ] \approx -\frac{1}{N} \sum_{i=1}^N \log q_K\left(x_i\right), \quad x_i \sim q_K(x)$,  and obtain the upper bound from $\tau$ in Kolmogorov-Smirnov distance.

\begin{algorithm}
    \caption{Steering Occupation Measure for a Large Group of Agents}
    \label{alg:3}
    \begin{algorithmic}[1]
        \Require The number of agents $N \in \mathbb{N}_{0}$; the maximal time step $K$; the system parameter $a(k)$ for $k = 0, \cdots, K-1$; the initial occupation measure $dq_{0}(x)$; the specified terminal occupation measure $d\tau(x)$.
        \Ensure Control inputs for the $i_\text{th}$ target $u_{i}(k)$, $k = 0, \cdots, K-1$, $i = 1, \cdots, N$.
        \State $k \Leftarrow 0$
    \While{$k < K$ and $e(k) \notin \mathbb{V}^{2n}_{++}$}
        \State Calculate $\mathscr{X}(k)$ using \eqref{momentsystem} if $k > 0$ or \eqref{XK} if $k = 0$
        \State Calculate $e(k)$ using \eqref{ek}
        \If{$e(k) \in \mathbb{V}^{2n}_{++}$}
        \State Calculate the states of the moment system $\mathscr{X}(i)$ for $i = k + 1, \cdots, K-1$ using \eqref{Xksum} with $\omega_{k} = \cdots = \omega_{K-1}$ \label{Step61}
        \State Calculate the controls of the moment system $\mathscr{U}(i)$ for $i = k, \cdots, K-1$ using \eqref{momentsystem}
        \If{$\exists i, \mathscr{U}(i) \notin \mathbb{V}^{2n}_{++}$}
        \State Back to Step \ref{Step6}, increase $\omega_{K-1}$ and determine other $\omega_{0}, \cdots, \omega_{K-2}$ by \eqref{Adjustk}
        \EndIf
        \State Optimize the cost function \eqref{Jr} and obtain the analytic estimates of the densities $\hat{p}_{i}(u)$ for $i = k, \cdots, K-1$
        \State  Sample the control inputs $u_{i}(j)$ from $\hat{p}_{i}(u)$ by acceptance-rejection sampling for all agents at time step $j = k, \cdots, K-1$ 
        \Else
            \State $u_{i}(k) = 0, i = 1, \cdots, N$
        \EndIf
        \State Calculate the power moments of the system state $x(k+1)$, i.e., $\mathscr{X}(k+1)$
        \State $k \Leftarrow k+1$
    \EndWhile
    \end{algorithmic}
\end{algorithm}

\section{Numerical examples}

In the previous sections, we proposed algorithms for steering a vast group of agents, either characterized as a probability density function or an occupation measure. In this section, we perform numerical simulations on different types of distributions supported on $\mathbb{R}$, to validate our proposed algorithms. We first simulate the steering of group agents as a probability density function. 

\subsection{Density steering of a vast group of agents}

In Example~1, we consider a four-step density steering task. The initial density is chosen as
\begin{equation}
    q_{0}(x) = \frac{1}{2} e^{\left( - x - 1 - e^{-(x+1)} \right)}
             + \frac{1}{2} e^{\left( - x + 1 - e^{-(x-1)} \right)},
    \label{q01}
\end{equation}
which is a mixture of two Gumbel densities. The terminal density is specified as
\begin{equation}
    \tau(x) = \frac{1}{2\sqrt{2\pi}}
              e^{\left( -\frac{(x-2)^2}{2} \right)}
            + \frac{1}{2\sqrt{2\pi}}
              e^{\left( -\frac{(x+2)^2}{2} \right)},
    \label{qt1}
\end{equation}
that is, a mixture of two Gaussian densities.

\begin{figure}[htbp]
\centering
\includegraphics[scale=0.3]{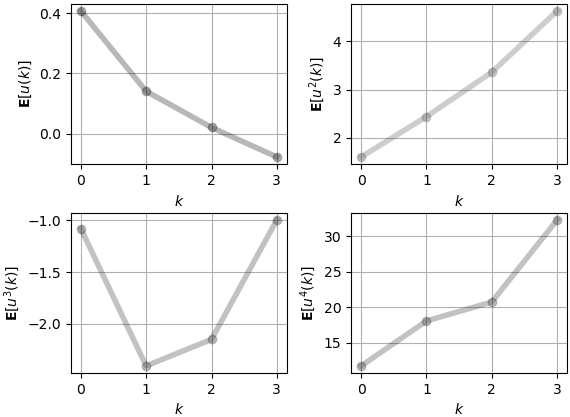}
\centering
\caption{$\mathscr{U}(k)$ at time steps $k = 0, 1, 2, 3$. The upper left figure shows $\mathbf{E}\left[ u(k)\right]$. The upper right one shows $\mathbf{E}\left[ u^{2}(k)\right]$. The lower left one shows $\mathbf{E}\left[ u^{3}(k)\right]$ and the lower right one shows $\mathbf{E}\left[ u^{4}(k)\right]$.}
\label{fig2}
\end{figure}

\begin{figure}[htbp]
\centering
\includegraphics[scale=0.3]{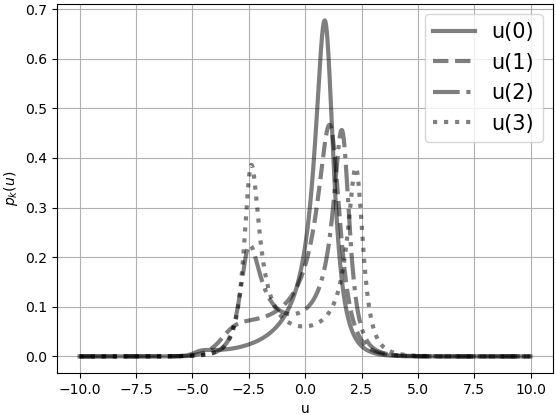}
\centering
\caption{Realized control inputs $u(k)$ by $\mathscr{U}(k)$ for $k = 0, 1, 2, 3$, which are obtained by our proposed control scheme.}
\label{fig3}
\end{figure}

The system parameters $a(k), k = 0, \cdots, 3$ are i.i.d. samples drawn from the uniform distribution $U[0.5, 0.7]$. The controls of the moment system, i.e., $\mathscr{U}(k)$ for $k = 0, 1, 2, 3$ are given in Figure \ref{fig2}. We note that by our proposed algorithm, $\mathscr{X}(k), \mathscr{U}(k) \in \mathbb{V}^{2n}_{++}$, which makes it feasible for us to realize the controls. The realized control inputs are given in Figure \ref{fig3}.

\begin{figure}[htbp]
\centering
\includegraphics[scale=0.3]{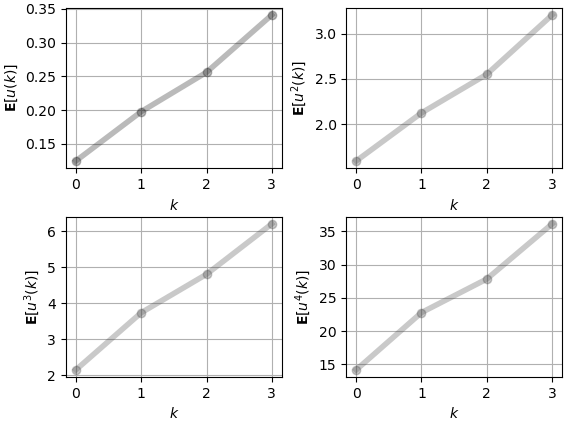}
\centering
\caption{$\mathscr{U}(k)$ at time steps $k = 0, 1, 2, 3$.}
\label{fig5}
\end{figure}

\begin{figure}[htbp]
\centering
\includegraphics[scale=0.3]{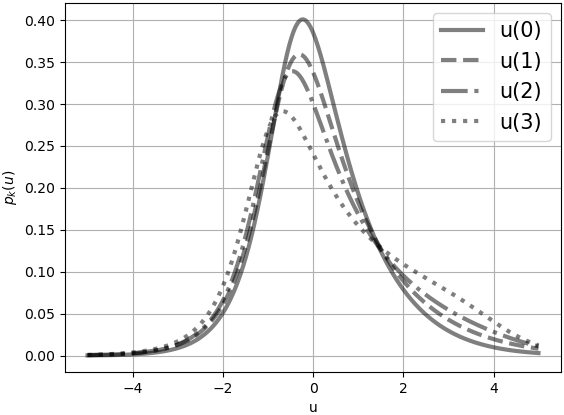}
\centering
\caption{Realized control inputs $u(k)$.}
\label{fig6}
\end{figure}

In Example 2, we simulate a steering problem in four steps where the initial density function is a Gaussian and the terminal density function is a mixture of generalized logistic distributions. The initial one is chosen as
$$
     q_{0}(x) = \frac{1}{\sqrt{2\pi}}e^{\frac{x ^{2}}{2}}.
$$
and the terminal one is specified as
\begin{equation}
     \tau(x) = \frac{0.4 \cdot 2e^{-x+1}}{(1 + e^{-x+1})^{3}} + \frac{0.6 \cdot 3e^{-x - 2}}{(1 + e^{-x - 2})^{4}}. 
\label{qt2}
\end{equation}
The system parameters $a(k), k = 0, \cdots, 3$ are i.i.d. samples drawn from the uniform distribution $U[0.5, 0.7]$. The controls of the moment system, i.e., $\mathscr{U}(k)$ for $k = 0, 1, 2, 3$ are given in Figure \ref{fig5}. The realized controls are shown in Figure \ref{fig6}.

\subsection{Occupation measure steering of a vast group of agents}

We simulate examples on occupation measure steering in this part of section. In Example 3, we steer $2000$ agents to a specified occupation measure. The initial states of each agent $x_{i}$ is drawn i.i.d. from the initial density \eqref{q01}. The specified terminal occupation measure consists of $2000$ i.i.d. samples drawn from the terminal distribution \eqref{qt1}. The system parameter $a(k), k = 0, \cdots, 3$ are i.i.d. samples drawn from the uniform distribution $U[0.5, 0.7]$. 

\begin{figure}[htbp]
\centering
\includegraphics[scale=0.35]{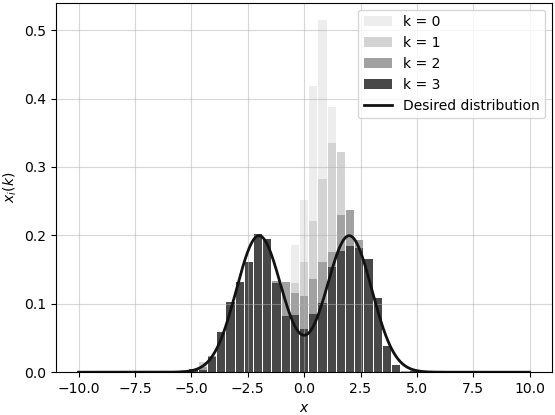}
\centering
\caption{The histograms of the system states $x(k)$ of 2000 Monte-Carlo simulations at time step $k = 0, 1, 2, 3$, together with the desired terminal distribution \eqref{qt1}. We note that the histogram of system state at step $k=3$ is very close to \eqref{qt1}.}
\label{fig14}
\end{figure}

Figure \ref{fig14} illustrates the histograms of the occupation measure for the states of the $2000$ agents at time steps $k = 0, 1, 2, 3$. The power moments of order $1$ to $4$ for the terminal occupation measure, obtained using our proposed algorithm, are $-1.27\times 10^{-2}, 4.87, -4.61\times 10^{-1}, 43.38$ respectively. Notably, these values are close to the desired terminal distribution, with power moments of order $1$ to $4$ being $0, 5, 0, 43$ respectively. 

In Example 4, we steer $2000$ agents to a specified occupation measure. The initial states of each agent $x_{i}$ is drawn i.i.d. from the Gaussian distribution \eqref{q01}. The specified terminal occupation measure consists of $2000$ i.i.d. samples drawn from the terminal distribution \eqref{qt2}. The system parameter $a(k), k = 0, \cdots, 3$ are i.i.d. samples drawn from the uniform distribution $U[0.5, 0.7]$. 

\begin{figure}[htbp]
\centering
\includegraphics[scale=0.35]{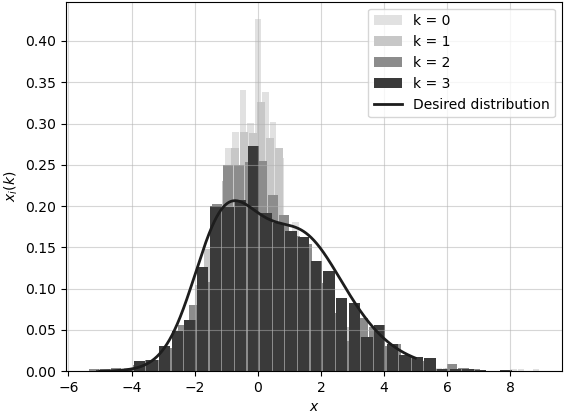}
\centering
\caption{The histograms of the system states $x(k)$ of 2000 Monte-Carlo simulations at time step $k = 0, 1, 2, 3$, together with the desired terminal distribution \eqref{qt2}.}
\label{fig16}
\end{figure}

Figure \ref{fig16} gives the histograms of the occupation measure of the states of the $2000$ agents at $k = 0, 1, 2, 3$. We note that the power moments of order $1$ to $4$ of \eqref{qt2} are $0.5, 3.88, 8.8, 52.8$ respectively. The power moments of order $1$ to $4$ of the terminal occupation measure by our proposed algorithm are $0.54, 3.97, 9.0, 52.79$ respectively, which are quite close to the specified ones. 

\section{Concluding remarks}

In this paper, we propose to use moments for the problem of steering a vast group of agents. The vast group of agents are characterized as probability density functions and occupation measures. We first treat the density steering problem. Without assuming the initial and terminal density to fall within specific function classes, the original problem is infinite-dimensional and intractable. As to treat this problem, we propose a moment system representation of the original system, and reduce the original problem to the control of the moment system. Different from the conventional control problems, the elements of the control inputs are in the system matrix, and we have to ensure the Hankel matrix of the control inputs at each time step to be positive definite. Since it is not treatable by the existing control schemes including the commonly used optimal control, we propose an empirical control scheme to treat this problem. By doing this, it is feasible for us to realize the control inputs for the original system. The previous results were presented at IFAC World Congress 2023 \cite{wu2023density}. In this paper, we consider more fundamental issues of the group steering problems and propose algorithms for real-world group steering problems. We propose that with the number of moment terms used approaching infinity, the terminal distribution obtained by the proposed scheme is almost everywhere equal to the desired one. We also propose an error upper bound of the terminal density by using our proposed algorithm, in the sense of the Kolmogorov-Smirnov distance. Based on the proposed density-steering algorithm, we put forward an algorithm steering an arbitrary occupation measure representing the vast group of agents to another arbitrary one. 
\bibliographystyle{plain}
\bibliography{autosam}

\end{document}